# The Approximate Solution of Newell-Whitehead-Segel and Fisher Equations Using The Adomian Decomposition Method

Asmaa A. Aswhad and Aqeel Falih Jaddoa
Department of Mathematics, College of Eduction for Pure Science/Ibn Al- Haitham, University of Baghdad.

الخلاصة

أستعملنا في البحث المقدم طريقة أدومين التحليلية لأيجاد الحل التقريبي لبعض حالات معادلة نيول- وايتهيد- سيكل التفاضلية غير الخطية والتي وجد الحل المضبوط لها سابقا بأستعمال طريقة هوموتوبي المتخلخلة والطريقة التكرارية ثم قارنا بين النتائج.

الكلمات المفتاحية: طريقة أدومين التحليلية، معادلة نيول- وايتهيد- سيكل، طريقة هوموتوبي المتخلخلة،طريقة تكرارية.

## ABSTRACT

In the present work, we use the Adomian Decomposition method to find the approximate solution for some cases of the Newell whitehead-segel nonlinear differential equation which was solved previously with exact solution by the Homotopy perturbation and the Iteration methods, then we compared the results.

## INTRODUCTION

The topic of the Adomian decomposition method has been rapidly growing in recent years. The concept of this method was first introduced by G. Adomian at the beginning of 1980's [1,2]. An advantage of the decomposition method is that it can provide analytical approximation to a rather wide class nonlinear and Stochastic equations without linearization, Perturbation, closure approximations or discretization methods which can result in massive numerical computation, [3]. Some of researchers modified Adomian decomposition method to solve generalize fifth order Kdv equations, [4]. This paper presents the use of the Adomian decomposition method for solving nonlinear partial differential equations analytically with initial conditions and compare the result with other method. We solved the Fisher equation and Newell- whitehead –segel equation by using Adomian decomposition method and compare it with the exact solution of the Iteration and the Homotopy Perturbation methods, [5, 6]. Some examples are prepared to illustrate these considerations.

**2. Adomian decomposition method,[3].**

We are with an equation $Fu(t)=g(t)$, where F represents a general nonlinear ordinary differential operator involving both linear and nonlinear terms.

The linear terms is decomposed into L+R, where L is an easily invertible and R is the remainder of a linear operator.

For convenience, L may be taken as the highest order derivative. Thus, the equation could be written as:

$$Lu + Ru + Nu = g, \quad \ldots(1)$$

Where Nu represents the nonlinear terms. Solving for Lu,

$$Lu = g - Ru - Nu. \quad \ldots(2)$$

Because L is invertible, an equivalent expression is

$$L^{-1}Lu = L^{-1}g - L^{-1}Ru - L^{-1}Nu. \quad \ldots(3)$$

If this corresponds to an initial –value problem, the integral operator $L^{-1}$ may be regarded as definite integrals from $t_0$ to t. If L is a second- order operator, $L^{-1}$ is twofold integration operator and

$$L^{-1}Lu = u - u(t_0) - (t - t_0)\acute{u}(t_0).$$

For boundary value problems indefinite integrations are used and the constants are evaluated from the given conditions.

Solving (3) for u yields

$$u = A + Bt + L^{-1}g - L^{-1}Ru - L^{-1}Nu. \quad \ldots(4)$$

The nonlinear term Nu will be equated to $\sum_{n=0}^{\infty} A_n$, where the $A_n$ are special polynomials defined as:

$$A_0 = f(u_0)$$

$$A_1 = u_1 \left(\frac{d}{du_0}\right) f(u_0)$$

$$A_2 = u_2 \left(\frac{d}{du_0}\right) f(u_0) + \left(\frac{u_1^2}{2}\right)\left(\frac{d^2}{du_0^2}\right) f(u_0)$$

$$A_3 = u_3 \left(\frac{d}{du_0}\right) f(u_0) + u_1 u_2 \left(\frac{d^2}{du_0^2}\right) f(u_0)$$

$$+ \left(\frac{u_1^3}{3}\right)\left(\frac{d^3}{du_0^3}\right) f(u_0)$$

.
.
.

Alternative definitions and formulas have been discussed in [7], and elsewhere one form of An is

$$A_n = \left(\frac{1}{n}\right) \sum_{v=1}^{n} c(v,n) \frac{d^v f}{du^v}$$

Where the second index in the coefficient is the order of the derivative and the first index progresses from 1 to n along with order of derivative.

Now u will be decomposed into $\sum_{n=0}^{\infty} u_n$ with $u_0$ identified as $A + Bt + L^{-1}g$ :

$$\sum_{n=0}^{\infty} u_n = u_0 - L^{-1}R \sum_{n=0}^{\infty} u_n - L^{-1} \sum_{n=0}^{\infty} A_n.$$

Consequently, we can write:

$$u_1 = -L^{-1}Ru_0 - L^{-1}A_0$$
$$u_2 = -L^{-1}Ru_1 - L^{-1}A_1 \quad (5)$$

.

$$\dot{u}_{n+1} = L^{-1}Ru_n - L^{-1}A_n$$

The polynomials $A_n$ are generated for nonlinearty so that $A_0$ depends only on $u_0$, $A_1$ depends only on $u_0$ and $u_1$, $A_2$ calcuable, and $u = \sum_{n=0}^{\infty} u_n$.

## 3. Newell-Whitehead- Segel equation

The Newell-Whitehead-Ssegel equation is written as:
$$u_t = ku_{xx} + au - bu^q \quad \ldots(6)$$
Where a, b and k are real numbers with $k > 0$, and q is a positive integer. To illustrate the Adomian method on this equation, three cases of nonlinear diffusion equation are presented

**Case 1. The Fisher's equation** (a=1, b=1, k=1 and q=2)

The Fisher's equation is
$$u_t = u_{xx} + u(1-u), \quad \ldots(7)$$
With a constant initial condition
$$u(x,0) = \alpha = u_0 \quad \ldots(8)$$
then $u_t - u_{xx} - u(1-u) = 0$, we obtain

$$\sum_{n=0}^{\infty} u_t = u_0 - L^{-1}R \sum_{n=0}^{\infty} u_n - L^{-1} \sum_{n=0}^{\infty} A_n$$

where $A_0 = f(u_0) = -u_0(1 - u_0)$
and from (5)
$$u_1 = -L^{-1}Ru_0 - L^{-1}A_0$$
$$= -\int_0^t Ru_0\, dt - \int_0^t A_0 dt = -\int_0^t -u_{0_{xx}}\, dt + \int_0^t u_0(1-u_0)dt$$
$$= 0 + u_0(1-u_0)t = \alpha(1-\alpha)t$$
$$A_1 = u_1\left(\frac{d}{du_0}\right)f(u_0)$$
$$= \alpha(1-\alpha)t(-1+2\alpha) = -\alpha(1-\alpha)(1-2\alpha)t$$
$$u_2 = -\int_0^t Ru_1 dt - \int_0^t A_1\, dt$$
$$= -\int_0^t -u_{1_{xx}}\, dt + \int_0^t \alpha(1-\alpha)(1-2\alpha)t = 0 + \alpha(1-\alpha)(1-2\alpha)\frac{t^2}{2!}$$
$$A_2 = u_2 \frac{d}{du_0}f(u_0) + \left(\frac{u_1^2}{2}\right)\frac{d^2}{du_0^2}f(u_0)$$
$$= -\alpha(1-\alpha)(-1+2\alpha)^2 \frac{t^2}{2} + \frac{2\alpha^2(1-\alpha)^2 t^2}{2}$$
$$= \alpha(1-\alpha)[-(1-2\alpha)^2 + 2\alpha(1-\alpha)]\,\frac{t^2}{2!}$$
$$= -\alpha(1-\alpha)[1 - 6\alpha + 6\alpha^2]\,\frac{t^2}{2!}$$

$$u_3 = -\int_0^t Ru_{2_{xx}} dt - \int_0^t A_2\, dt$$

$$= -\int_0^t -u_{2_{xx}}\, dt - \int_0^t -\alpha(1-\alpha)[1 - 6\alpha + 6\alpha^2]\,\frac{t^2}{2!}\,dt$$

$$= \alpha(1-\alpha)(1 - 6\alpha + 6\alpha^2)\frac{t^3}{3!}$$

$$A_3 = u_3 \frac{d}{du_0} f(u_0) + u_1 u_2 \frac{d^2}{du_0^2} f(u_0) + \left(\frac{u_1^3}{3!}\right)\frac{d^3}{du_0^3} f(u_0)$$

$$= -\alpha(1-\alpha)(1 - 6\alpha + 6\alpha^2)\frac{t^3}{3!}(1 - 2\alpha) + 2\alpha(1-\alpha)t\alpha(1-\alpha)(1 - 2\alpha)\frac{t^2}{2}$$

$$= -\alpha(1-\alpha)(1 - 2\alpha)(1 - 6\alpha + 6\alpha^2)\frac{t^3}{3!} + 6\alpha^2(1-\alpha)^2(1 - 2\alpha)\frac{t^3}{6}$$

$$= -\alpha(1-\alpha)(1 - 2\alpha)[1 - 12\alpha + 12\alpha^2]\,\frac{t^3}{3!}$$

$$u_4 = -\int_0^t Ru_{3_{xx}} dt - \int_0^t A_3\, dt$$

$$= 0 - \int_0^t -\alpha(1-\alpha)(1 - 2\alpha)[1 - 12\alpha + 12\alpha^2]\,\frac{t^3}{3!}\,dt$$

$$= \alpha(1-\alpha)(1 - 2\alpha)(1 - 12\alpha + 12\alpha^2)\frac{t^4}{4!}$$

Since, u(x, t)=$\sum_{n=0}^{\infty} u_n = u_0 + u_1 + u_2 + \cdots$

Then $u(x,t) = \alpha(1-\alpha)t + \alpha(1-\alpha)(1-2\alpha)\frac{t^2}{2!} + \alpha(1-\alpha)(1 - 6\alpha + 6\alpha^2)\frac{t^3}{3!} + \cdots$

**Case 2: In equation (6) for a=1, b=1, k=1, q=2. The Newell- whitehead segel equation is written as:**

$$u_t = u_{xx} + u - u^2 \qquad \ldots(10)$$

With the initial condition

$$u(x, 0) = \frac{1}{\left(1+e^{x/\sqrt{6}}\right)^2} = u_0 \qquad \ldots (11)$$

Then $u_t - u_{xx} - (u - u^2) = 0$,

$$\sum_{n=0}^{\infty} u_n = u_0 - L^{-1} R \sum_{n=0}^{\infty} u_n - L^{-1} \sum_{n=0}^{\infty} A_n$$

$$A_0 = f(u_0) = Nu_0 = -(u_0 - u_0^2)$$

$$= -\left(\frac{1}{\left(1+e^{x/\sqrt{6}}\right)^2} - \frac{1}{\left(1+e^{x/\sqrt{6}}\right)^4}\right)$$

$$u_1 = -L^{-1}Ru_0 - L^{-1}A_0 = -\int_0^t -u_{0_{xx}}\, dt + \int_0^t A_0 dt$$

$$= -\int_0^t \frac{e^{x/\sqrt{6}}\left(2e^{x/\sqrt{6}} - 1\right)}{3\left(1 + e^{x/\sqrt{6}}\right)^4}\, dt + \int_0^t \left(\frac{1}{\left(1 + e^{x/\sqrt{6}}\right)^2} - \frac{1}{\left(1 + e^{x/\sqrt{6}}\right)^4}\right) dt$$

$$= \left[\frac{e^{x/\sqrt{6}}\left(2e^{x/\sqrt{6}} - 1\right)}{3\left(1 + e^{x/\sqrt{6}}\right)^4} + \frac{1}{\left(1 + e^{x/\sqrt{6}}\right)^2} - \frac{1}{\left(1 + e^{x/\sqrt{6}}\right)^4}\right] t = \frac{5e^{x/\sqrt{6}}}{3\left(1 + e^{x/\sqrt{6}}\right)^3} t$$

$$A_1 = u_1 \frac{df(u_0)}{du_0}$$

$$= -\frac{5e^{x/\sqrt{6}}}{3\left(1 + e^{x/\sqrt{6}}\right)^3}\left(1 - \frac{2}{\left(1 + e^{x/\sqrt{6}}\right)^2}\right) t$$

$$u_{1x} = \frac{5e^{x/\sqrt{6}}\left(1 - 2e^{x/\sqrt{6}}\right)}{3\sqrt{6}\left(1 + e^{x/\sqrt{6}}\right)^4} t$$

$$u_{1xx} = \frac{5e^{x/\sqrt{6}}\left(1 - 7e^{x/\sqrt{6}} + 4\left(e^{x/\sqrt{6}}\right)^2\right)}{18\left(1 + e^{x/\sqrt{6}}\right)^5} t$$

$$u_2 = -L^{-1}Ru_1 - L^{-1}A_1 = -\int_0^t -u_{1_{xx}}\, dt + \int_0^t A_1 dt$$

$$= \int_0^t \frac{5e^{x/\sqrt{6}}\left(1 - 7e^{x/\sqrt{6}} + 4\left(e^{x/\sqrt{6}}\right)^2\right)}{18\left(1 + e^{x/\sqrt{6}}\right)^5} t\, dt$$

$$+ \int_0^t \frac{5e^{x/\sqrt{6}}\left(\left(1 + e^{x/\sqrt{6}}\right)^2 - 2\right)}{3\left(1 + e^{x/\sqrt{6}}\right)^5} t\, dt$$

$$= \frac{5e^{x/\sqrt{6}}\left(1 - 7e^{x/\sqrt{6}} + 4\left(e^{x/\sqrt{6}}\right)^2\right) + 6\left(-1 + 2e^{x/\sqrt{6}} + \left(e^{x/\sqrt{6}}\right)^2\right)}{18\left(1 + e^{x/\sqrt{6}}\right)^5} \frac{t^2}{2}$$

$$= \frac{25e^{x/\sqrt{6}}\left(-1 + e^{x/\sqrt{6}} + 2e^{2x/\sqrt{6}}\right)}{18\left(1 + e^{x/\sqrt{6}}\right)^5} \frac{t^2}{2}$$

$$= \frac{25e^{x/\sqrt{6}}\left(\left(-1+2e^{x/\sqrt{6}}\right)\left(1+e^{x/\sqrt{6}}\right)\right)}{18\left(1+e^{x/\sqrt{6}}\right)^5}\frac{t^2}{2} = \frac{25e^{x/\sqrt{6}}\left(-1+2e^{x/\sqrt{6}}\right)}{18\left(1+e^{x/\sqrt{6}}\right)^4}\frac{t^2}{2}$$

$$A_2 = u_2 \frac{df(u_0)}{du_0} + \left(\frac{u_1^2}{2}\right)\frac{d^2f(u_0)}{du_0^2}$$

$$= \left[\frac{25e^{x/\sqrt{6}}\left(-1+2e^{x/\sqrt{6}}\right)\left(2-\left(1+e^{x/\sqrt{6}}\right)^2\right)}{18\left(1+e^{x/\sqrt{6}}\right)^6} + \frac{100\left(e^{x/\sqrt{6}}\right)^2}{18\left(1+e^{x/\sqrt{6}}\right)^6}\right]\frac{t^2}{2}$$

$$= \frac{25e^{x/\sqrt{6}}\left[\left(-1+2e^{x/\sqrt{6}}\right)\left(1-2e^{x/\sqrt{6}}-\left(e^{x/\sqrt{6}}\right)^2\right)+4e^{x/\sqrt{6}}\right]}{18\left(1+e^{x/\sqrt{6}}\right)^6}\frac{t^2}{2}$$

$$= \frac{25e^{x/\sqrt{6}}\left(-1+8e^{x/\sqrt{6}}-3\left(e^{x/\sqrt{6}}\right)^2-2\left(e^{x/\sqrt{6}}\right)^3\right)}{18\left(1+e^{x/\sqrt{6}}\right)^6}\frac{t^2}{2}$$

$$u_{2x} = \frac{25\left(\left(e^{x/\sqrt{6}}\right)^3+7\left(e^{x/\sqrt{6}}\right)^2-e^{x/\sqrt{6}}\right)}{3\sqrt{6}\left(1+e^{x/\sqrt{6}}\right)^5}$$

$$u_{2xx} = \frac{25\left(8\left(e^{x/\sqrt{6}}\right)^4 - 33\left(e^{x/\sqrt{6}}\right)^3 - 18\left(e^{x/\sqrt{6}}\right)^2 - e^{x/\sqrt{6}}\right)}{108\left(1+e^{x/\sqrt{6}}\right)^6}\frac{t^2}{2}$$

$$u_3 = -L^{-1}Ru_2 - L^{-1}A_2 = -\int_0^t -u_{2xx}\,dt - \int_0^t A_2\,dt$$

$$= \frac{25e^{x/\sqrt{6}}\left(8\left(e^{x/\sqrt{6}}\right)^3 - 33\left(e^{x/\sqrt{6}}\right)^2 + 18e^{x/\sqrt{6}} - 1 + 6 - 48e^{x/\sqrt{6}} + 18\left(e^{x/\sqrt{6}}\right)^2 + 12\left(e^{x/\sqrt{6}}\right)\right)}{216\left(1+e^{x/\sqrt{6}}\right)^6}$$

$$= \frac{25e^{x/\sqrt{6}}\left(20\left(e^{x/\sqrt{6}}\right)^3 - 15\left(e^{x/\sqrt{6}}\right)^2 - 30e^{x/\sqrt{6}} + 5\right)}{216\left(1+e^{x/\sqrt{6}}\right)^6}\frac{t^3}{3}$$

$$= \frac{125e^{x/\sqrt{6}}\left(1+e^{x/\sqrt{6}}\right)\left(4\left(e^{x/\sqrt{6}}\right)^2 - 7e^{x/\sqrt{6}} + 1\right)}{216\left(1+e^{x/\sqrt{6}}\right)^6}\frac{t^3}{3}$$

$$= \frac{125e^{x/\sqrt{6}}\left(4\left(e^{x/\sqrt{6}}\right)^2 - 7e^{x/\sqrt{6}} + 1\right)}{216\left(1+e^{x/\sqrt{6}}\right)^5}\frac{t^3}{3}$$

Since $u(x, t) = \sum_{n=0}^{\infty} u_n = u_0 + u_1 + u_2 + \cdots$ , then

$$u(x,t) = \frac{1}{\left(1+e^{x/\sqrt{6}}\right)^2} + \frac{5e^{x/\sqrt{6}}}{3\left(1+e^{x/\sqrt{6}}\right)^3}t + \frac{25}{18}\left(\frac{e^{x/\sqrt{6}}\left(-1+2e^{x/\sqrt{6}}\right)}{\left(1+e^{x/\sqrt{6}}\right)^4}\right)\frac{t^2}{2} + \cdots \quad \ldots(12)$$

**Case 3: In this case we will examine the Newell- whitehead segel equation for a=3, b=4, k=1, q=3**

$$u_t = u_{xx} + 3u - 4u^3 \quad \ldots(13)$$

With the initial condition

$$u_0 = u(x, 0) = \sqrt{\frac{3}{4}} \; \frac{e^{\sqrt{6}x}}{e^{\sqrt{6}x}+e^{\frac{\sqrt{6}}{2}x}} \quad \ldots(14)$$

$u_t - u_{xx} - 3u + 4u^3 = 0$

From $A_0 = f(u_0) = -3u_0 + 4u_0^{\,3}$

$$u_{0_x} = \sqrt{\frac{18}{4}} \; \frac{\frac{1}{2}e^{\frac{3\sqrt{6}}{2}x}}{\left(e^{\sqrt{6}x} + e^{\frac{\sqrt{6}}{2}x}\right)^2}$$

$$u_{0_{xx}} = \frac{6\sqrt{3}}{4\sqrt{4}} \; \frac{e^{\frac{3\sqrt{6}}{2}x}\left(-e^{\sqrt{6}x} + e^{\frac{\sqrt{6}}{2}x}\right)}{\left(e^{\sqrt{6}x} + e^{\frac{\sqrt{6}}{2}x}\right)^3}$$

We obtain

$$u_1 = -\int_0^t -u_{0_{xx}}\, dt - \int_0^t (-3u_0 + 4u_0^{\,3})dt$$

$$= -\int_0^t -u_{0_{xx}}\, dt + \int_0^t (3u_0 - 4u_0^{\,3})dt$$

$$= \int_0^t \frac{6\sqrt{3}}{4\sqrt{4}} \; \frac{e^{\frac{3\sqrt{6}}{2}x}\left(-e^{\sqrt{6}x} + e^{\frac{\sqrt{6}}{2}x}\right)}{\left(e^{\sqrt{6}x} + e^{\frac{\sqrt{6}}{2}x}\right)^3}\, dt$$

$$+ \int_0^t \left[ 3\sqrt{\frac{3}{4}} \; \frac{e^{\sqrt{6}x}}{e^{\sqrt{6}x}+e^{\frac{\sqrt{6}}{2}x}} - 4\left(\sqrt{\frac{3}{4}}\right)^3 \frac{\left(e^{\sqrt{6}x}\right)^3}{\left(e^{\sqrt{6}x}+e^{\frac{\sqrt{6}}{2}x}\right)^3} \right] dt$$

$$= \left[ \frac{3}{2}\sqrt{\frac{3}{4}} \frac{e^{\frac{3\sqrt{6}}{2}x}(-e^{\sqrt{6}x} + e^{\frac{\sqrt{6}}{2}x})}{\left(e^{\sqrt{6}x} + e^{\frac{\sqrt{6}}{2}x}\right)^3} + 3\sqrt{\frac{3}{4}} \frac{e^{\sqrt{6}x}}{e^{\sqrt{6}x} + e^{\frac{\sqrt{6}}{2}x}} \right.$$

$$\left. - 3\sqrt{\frac{3}{4}} \frac{e^{3\sqrt{6}x}}{\left(e^{\sqrt{6}x} + e^{\frac{\sqrt{6}}{2}x}\right)^3} \right] t$$

$$= \frac{3}{2}\sqrt{\frac{3}{4}} \left[ \frac{e^{\frac{3\sqrt{6}}{2}x}\left(-e^{\sqrt{6}x} + e^{\frac{\sqrt{6}}{2}x}\right) + 2 - e^{\sqrt{6}x}(e^{\sqrt{6}x} + e^{\frac{\sqrt{6}}{2}x})^2 - 2e^{3\sqrt{6}x}}{\left(e^{\sqrt{6}x} + e^{\frac{\sqrt{6}}{2}x}\right)^3} \right] t$$

$$= \frac{3}{2}\sqrt{\frac{3}{4}} \left[ \frac{-e^{\frac{5\sqrt{6}}{2}x} + e^{2\sqrt{6}x} + 2e^{3\sqrt{6}x} + 4e^{\frac{5\sqrt{6}}{2}x} + 2e^{2\sqrt{6}x} - 2e^{3\sqrt{6}x}}{\left(e^{\sqrt{6}x} + e^{\frac{\sqrt{6}}{2}x}\right)^3} \right] t$$

$$= \frac{3}{2}\sqrt{\frac{3}{4}} \frac{3e^{\frac{5\sqrt{6}}{2}x} + 3e^{2\sqrt{6}x}}{\left(e^{\sqrt{6}x} + e^{\frac{\sqrt{6}}{2}x}\right)^3} t$$

$$u_1(x,t) = \frac{9}{2}\sqrt{\frac{3}{4}} \frac{e^{\sqrt{6}x} e^{\frac{\sqrt{6}}{2}x}}{\left(e^{\sqrt{6}x} + e^{\frac{\sqrt{6}}{2}x}\right)^2} t$$

$$A_1 = u_1 \frac{df(u_0)}{du_0} = u_1(-3 + 12u_0^2)$$

$$= \frac{9}{2}\sqrt{\frac{3}{4}} \frac{e^{\sqrt{6}x} e^{\frac{\sqrt{6}}{2}x}}{\left(e^{\sqrt{6}x} + e^{\frac{\sqrt{6}}{2}x}\right)^2} \left( -3 + 9 \frac{e^{2\sqrt{6}x}}{\left(e^{\sqrt{6}x} + e^{\frac{\sqrt{6}}{2}x}\right)^2} \right) t$$

$$= \frac{9}{2}\sqrt{\frac{3}{4}} \frac{e^{\sqrt{6}x} e^{\frac{\sqrt{6}}{2}x} \left( -3\left(e^{\sqrt{6}x} + e^{\frac{\sqrt{6}}{2}x}\right)^2 + 9e^{2\sqrt{6}x} \right)}{\left(e^{\sqrt{6}x} + e^{\frac{\sqrt{6}}{2}x}\right)^4} t$$

$$= \frac{9}{2}\sqrt{\frac{3}{4}} \ \frac{e^{\sqrt{6}x} e^{\frac{\sqrt{6}}{2}x}\left(6e^{2\sqrt{6}x} - 6e^{\frac{3\sqrt{6}}{2}x} - 3e^{\sqrt{6}x}\right)}{\left(e^{\sqrt{6}x} + e^{\frac{\sqrt{6}}{2}x}\right)^4} \ t$$

$$= \frac{27\sqrt{3}}{4} \ \frac{e^{\frac{3\sqrt{6}}{2}x}\left(2e^{2\sqrt{6}x} - 2e^{\frac{3\sqrt{6}}{2}x} - e^{\sqrt{6}x}\right)}{\left(e^{\sqrt{6}x} + e^{\frac{\sqrt{6}}{2}x}\right)^4} \ t$$

$$= \frac{27\sqrt{3}}{4} \ \frac{\left(2e^{\frac{7\sqrt{6}}{2}x} - 2e^{3\sqrt{6}x} - e^{\frac{5\sqrt{6}}{2}x}\right)}{\left(e^{\sqrt{6}x} + e^{\frac{\sqrt{6}}{2}x}\right)^4} \ t$$

from $u_1(x,t)$ on the previous page, we get

$$u_{1x} = \frac{9\sqrt{18}}{8} \ \frac{\left(-e^{\frac{5\sqrt{6}}{2}x} + e^{2\sqrt{6}x}\right)}{\left(e^{\sqrt{6}x} + e^{\frac{\sqrt{6}}{2}x}\right)^3}$$

$$u_{1xx} = \frac{27\sqrt{3}}{8} \ \frac{\left(e^{\frac{7\sqrt{6}}{2}x} - 4e^{3\sqrt{6}x} + e^{\frac{5\sqrt{6}}{2}x}\right)}{\left(e^{\sqrt{6}x} + e^{\frac{\sqrt{6}}{2}x}\right)^4}$$

$$u_2 = -\int_0^t -u_{1xx}\, dt - \int_0^t A_1 dt$$

$$= \frac{27\sqrt{3}}{8} \ \frac{\left(e^{\frac{7\sqrt{6}}{2}x} - 4e^{3\sqrt{6}x} + e^{\frac{5\sqrt{6}}{2}x}\right)}{\left(e^{\sqrt{6}x} + e^{\frac{\sqrt{6}}{2}x}\right)^4} \ \frac{t^2}{2}$$

$$- \frac{27\sqrt{3}}{4} \ \frac{\left(2e^{\frac{7\sqrt{6}}{2}x} - 2e^{3\sqrt{6}x} - e^{\frac{5\sqrt{6}}{2}x}\right)}{\left(e^{\sqrt{6}x} + e^{\frac{\sqrt{6}}{2}x}\right)^4} \ \frac{t^2}{2}$$

$$= \frac{27\sqrt{3}}{8} \frac{\left(-3e^{\frac{7\sqrt{6}}{2}x} + 3e^{\frac{5\sqrt{6}}{2}x}\right)}{\left(e^{\sqrt{6}x} + e^{\frac{\sqrt{6}}{2}x}\right)^4} \frac{t^2}{2}$$

$$= \frac{81}{4}\sqrt{\frac{3}{4}} \frac{\left(-e^{\frac{\sqrt{6}}{2}x}e^{3\sqrt{6}x} + e^{\frac{\sqrt{6}}{2}x}e^{2\sqrt{6}x}\right)}{\left(e^{\sqrt{6}x} + e^{\frac{\sqrt{6}}{2}x}\right)^4} \frac{t^2}{2}$$

$$= \frac{81}{4}\sqrt{\frac{3}{4}} \frac{e^{\frac{\sqrt{6}}{2}x}e^{\sqrt{6}x}\left(-e^{2\sqrt{6}x} + e^{\sqrt{6}x}\right)}{\left(e^{\sqrt{6}x} + e^{\frac{\sqrt{6}}{2}x}\right)^4} \frac{t^2}{2}$$

$$= \frac{81}{4}\sqrt{\frac{3}{4}} \frac{e^{\frac{\sqrt{6}}{2}x}e^{\sqrt{6}x}\left(e^{\sqrt{6}x} + e^{\frac{\sqrt{6}}{2}x}\right)\left(-e^{\sqrt{6}x} + e^{\frac{\sqrt{6}}{2}x}\right)}{\left(e^{\sqrt{6}x} + e^{\frac{\sqrt{6}}{2}x}\right)^4} \frac{t^2}{2}$$

$$= \frac{81}{4}\sqrt{\frac{3}{4}} \frac{e^{\frac{\sqrt{6}}{2}x}e^{\sqrt{6}x}\left(-e^{\sqrt{6}x} + e^{\frac{\sqrt{6}}{2}x}\right)}{\left(e^{\sqrt{6}x} + e^{\frac{\sqrt{6}}{2}x}\right)^3} \frac{t^2}{2}$$

$$A_2 = u_2 \frac{df(u_0)}{du_0} + \left(\frac{u_1^2}{2}\right)\frac{d^2 f(u_0)}{du_0^2}$$

$$= \frac{81}{4}\sqrt{\frac{3}{4}} \frac{e^{\frac{\sqrt{6}}{2}x}e^{\sqrt{6}x}\left(-e^{\sqrt{6}x} + e^{\frac{\sqrt{6}}{2}x}\right)}{\left(e^{\sqrt{6}x} + e^{\frac{\sqrt{6}}{2}x}\right)^3} \frac{t^2}{2}(-3 + 12u_0^2)$$

$$+ \frac{81}{4}\frac{3}{4} \frac{e^{2\sqrt{6}x}e^{\sqrt{6}x}}{\left(e^{\sqrt{6}x} + e^{\frac{\sqrt{6}}{2}x}\right)^4} \frac{t^2}{2} 24u_0$$

$$= \frac{243}{8}\sqrt{\frac{3}{4}} \frac{e^{3\sqrt{6}x}\left(-2e^{\frac{3\sqrt{6}}{2}x} + 10e^{\sqrt{6}x} - e^{\frac{\sqrt{6}}{2}x} - 1\right)}{\left(e^{\sqrt{6}x} + e^{\frac{\sqrt{6}}{2}x}\right)^5} t^2$$

$$u_3 = -\int_0^t -u_{2xx}\, dt - \int_0^t A_2 dt$$

By the same way we have

$$u_3 = \frac{243}{16}\sqrt{\frac{3}{4}}\frac{e^{\sqrt{6}x}e^{\frac{\sqrt{6}}{2}x}\left(-4e^{\sqrt{6}x}e^{\frac{\sqrt{6}}{2}x}+\left(e^{\sqrt{6}x}\right)^2+\left(e^{\frac{\sqrt{6}}{2}x}\right)^2\right)}{\left(e^{\sqrt{6}x}+e^{\frac{\sqrt{6}}{2}x}\right)^4}\frac{t^3}{3}$$

$u(x,t) = u_0 + u_1 + u_2 + u_3 + \cdots$

Then $u(x,t) = \sqrt{\frac{3}{4}}\frac{e^{\sqrt{6}x}}{e^{\sqrt{6}x}+e^{\frac{\sqrt{6}}{2}x}} + \frac{9}{2}\sqrt{\frac{3}{4}}\frac{e^{\sqrt{6}x}e^{\frac{\sqrt{6}}{2}x}}{\left(e^{\sqrt{6}x}+e^{\frac{\sqrt{6}}{2}x}\right)^2}t +$

$$\frac{81}{4}\sqrt{\frac{3}{4}}\frac{e^{\sqrt{6}x}e^{\frac{\sqrt{6}}{2}x}\left(-e^{\sqrt{6}x}+e^{\frac{\sqrt{6}}{2}x}\right)}{\left(e^{\sqrt{6}x}+e^{\frac{\sqrt{6}}{2}x}\right)^3}\frac{t^2}{2} + \quad \cdots \quad \ldots(15)$$

Now, from comparing our results and the results obtained in [5], [6] where the exact solution for Fisher's equation in [5], is

$$u(x,t) = \alpha + \alpha(1-\alpha)t + \alpha(1-\alpha)(1-2\alpha)\frac{t^2}{2!}$$
$$+ \alpha(1-\alpha)(1-6\alpha+6\alpha^2)\frac{t^3}{3!} + \cdots$$

that implies

$$u(x,t) = \frac{\alpha e^t}{1-\alpha+\alpha e^t} \quad \ldots(16)$$

And for Newell- Whitehead-Segel (case 2) in [6] is

$$u(x,t) = \frac{1}{\left(1+e^{\frac{x}{\sqrt{6}}}\right)^2} + \frac{5}{3}\frac{e^{\frac{x}{\sqrt{6}}}}{\left(1+e^{\frac{x}{\sqrt{6}}}\right)^3}t + \frac{25}{18}\left(\frac{e^{\frac{x}{\sqrt{6}}}\left(-1+2e^{\frac{x}{\sqrt{6}}}\right)}{\left(1+e^{\frac{x}{\sqrt{6}}}\right)^4}\right)\frac{t^2}{2} +$$

$$\frac{125}{216}\left(\frac{e^{\frac{x}{\sqrt{6}}}\left(4\left(e^{\frac{x}{\sqrt{6}}}\right)^2-7e^{\frac{x}{\sqrt{6}}}+1\right)}{\left(1+e^{\frac{x}{\sqrt{6}}}\right)^5}\right)\frac{t^3}{3} + \cdots$$

that equal

$$u(x,t) = \frac{1}{\left(1+e^{\frac{x}{\sqrt{6}}-\frac{5}{6}t}\right)^2} \quad \ldots(17)$$

For Newell-Whiteheat- Segel (case 3) in [6] is

$$u(x,t) = \sqrt{\frac{3}{4}} \frac{e^{\sqrt{6}x}}{e^{\sqrt{6}x} + e^{\frac{\sqrt{6}}{2}x}} + \frac{9}{2}\sqrt{\frac{3}{4}} \frac{e^{\sqrt{6}x}e^{\frac{\sqrt{6}}{2}x}}{\left(e^{\sqrt{6}x} + e^{\frac{\sqrt{6}}{2}x}\right)^2} t$$

$$+ \frac{81}{4}\sqrt{\frac{3}{4}} \frac{e^{\sqrt{6}x}e^{\frac{\sqrt{6}}{2}x}\left(-e^{\sqrt{6}x} + e^{\frac{\sqrt{6}}{2}x}\right)}{\left(e^{\sqrt{6}x} + e^{\frac{\sqrt{6}}{2}x}\right)^3} \frac{t^2}{2} + \cdots$$

that equal

$$u(x,t) = \sqrt{\frac{3}{4}} \frac{e^{\sqrt{6}x}}{e^{\sqrt{6}x} + e^{(\frac{\sqrt{6}}{2}x - \frac{9}{2}t)}} \qquad \ldots(18)$$

We note that the Adomian decomposition method gives an approximate solution which is very close to the exact solution.

## CONCLUSION

In the present work an approximate solution for Fisher׳s and Newell-Whitehead- Segel equations are obtained by Adomian decomposition method where the exact solution is known by Variational Iteration and the Homotopy perturbation methods, we concluded that the Adomian decomposition method one of the best methods to find the approximate solution because it gives a better results and closer to the exact solution as well as it's a simple method.